\documentclass[11pt,a4paper]{article}
\usepackage{amsmath,amssymb}
\usepackage{amsmath}

\usepackage{color}
\usepackage[colorlinks,linkcolor=blue,citecolor=red]{hyperref}

\DeclareMathSymbol{\bbbr}{\mathalpha}{AMSb}{"52}
\DeclareMathSymbol{\bbbc}{\mathalpha}{AMSb}{"52}

\newtheorem{theorem}{Theorem}

\newtheorem{corollary}[theorem]{Corollary}

\newtheorem{definition}[theorem]{Definition}

\newtheorem{lemma}[theorem]{Lemma}

\begin{document}

\title{ Duality for systems of conservation laws}

\author{{\Large Sergey I. Agafonov}\\
\\
Department of Mathematics,\\
S\~ao Paulo State University-UNESP,\\ S\~ao Jos\'e do Rio Preto, Brazil\\
e-mail: {\tt sergey.agafonov@gmail.com} }
\date{}
\maketitle
\unitlength=1mm

\vspace{1cm}

\begin{abstract} For one-dimensional systems of conservation laws admitting two additional conservation laws we assign a ruled hypersurface of codimension two in projective space. We call two such systems dual if the corresponding ruled hypersurfaces are dual. We show that a Hamiltonian system is autodual, its  ruled hypersurface sits in some quadric, and the generators of this ruled hypersurface form a Legendre submanifold
 with respect to the  contact structure on Fano variety of this quadric. We also give a complete geometric description of 3-component nondiagonalizable systems of Temple class: such systems admit two additional conservation laws, they are dual to systems with constant characteristic speeds, constructed via maximal rank 3-webs of curves in space.
\bigskip


\bigskip

\end{abstract}

%\vspace{-7mm}
%\newpage

%\tableofcontents

\section{Introduction}
The  most fundamental laws  of nature express the conservation of mass, momentum, energy, etc. The corresponding mathematical models, known as systems of conservation laws, were the object of extensive study from the classical period of mathematical physics, when Euler came up with equations for ideal fluid, to modern times (see \cite{J-75,RJ-83,S-99}).  Such systems revealed  a rich geometric structure
 (see for instance \cite{DN-83,T-85,S-94,AF-99}). This paper falls into the long-term trend of geometrization of mathematical physics.

In presence of 2-dimensional translation symmetry, systems of conservation laws reduce to a system of quasilinear partial differential equations (PDEs) in two independent variables
\begin{equation}
u^i_t=\frac{d}{dx}f^i(u)=v^i_j(u)  u^j_x, ~~~
v^i_j=\frac{\partial f^i}{\partial u^j},~~~i=1,...,n,
\label{CL}
\end{equation}
where $t$ is interpreted as the time and $x$ as the space variable. (In what follows, we sum by repeating indices.)

To the system of conservation laws (\ref{CL}), one associates a
 {\it congruence of lines} (see \cite{AF-96}), i.e. an $n$-parameter family of lines
\begin{equation}
Y^i=u^i\ Y^0+f^i(u)Y^{n+1}, ~~~ i=1,..., n
\label{congrn}
\end{equation}
in $(n+1)$-dimensional projective space $\mathbb P^{n+1}$ with homogeneous coordinates
$[Y^0:...:Y^{n+1}]$. This congruence is also viewed as an $n$-dimensional submanifold of Grassmannian $Gr(2,n+2)$.

The geometry behind this correspondence is the following. Any solution $u(x,t)$ to  (\ref{CL}), which is almost inevitably local for non-linear systems, defines a (germ of) surface $S$ in $(n+2)$-dimensional affine space $A^{n+2}$ with coordinates $(q^1,...,q^n,x,t)=(q,x,t)$, parametrized by $x,t$ via integrating the forms $dq^i:=u^i(x,t)dx+f^i(u(x,t))dt$ closed on the solutions. Taking the affine tangent space to this surface at a point $(q(x,t),x,t)\in S$, translating it to the origin, and projectivizing, we get the {\it Gauss map}
$$
G:S\to Gr(2,n+2).
$$
Now solutions to (\ref{CL}) are intrinsically interpreted as surfaces in $A^{n+2}$ whose Gauss map images lie in the Grassmannian submanifold (\ref{congrn}).

The correspondence allows one to translate   the basic concepts of the theory of systems
of conservation laws into geometrical language. In particular (see \cite{AF-96,AF-99}),
\begin{itemize}
\item rarefaction curves parametrize developable surfaces of the congruence;
\item shock curves describe ruled surfaces of the congruence, whose line generators meet a fixed congruence line;
\item reciprocal transformations manifest themselves as projective transformations of the congruence;
\item commuting systems correspond to congruences, whose focal submanifolds are related by a Peterson transformation;
\item linear degeneracy of a system means that the  developable surfaces of the congruence are conical and the codimension of the focal submanifold jumps up by one;
\item rarefaction curves are rectilinear (this property distinguishes systems of Temple class) if and only if the developable surfaces of the congruence are planar;
\item  a system possesses a Riemann invariant when the focal net of the corresponding focal submanifold is conjugate.
\end{itemize}

Phrased in this geometrical language, the classification of 2-component linearly degenerate systems and of systems of Temple class \cite{T-83} becomes an easy exercise in analytic geometry \cite{AF-96}.
The case of 3-component  systems  is
more involved. In presence of Riemann invariants, both linearly degenerate systems and systems of Temple class were described geometrically in \cite{AF-99}. The case of 3-component nondiagonalizable systems (i.e. not possessing a single Riemann invariant) revealed much stronger resistance. Three-component systems of Temple class were studied in  \cite{AF-99}:  the derived so-called structure equations and equations for characteristic speeds allowed to show that systems with the same structure equations are equivalent by reciprocal transformations; integrable and Hamiltonian systems were also identified . Integrability, linear degeneracy, and equivalence to some Hamiltonian system by reciprocal transformations were shown to be equivalent.

The congruence of a linearly degenerate 3-component system of Temple class turns out to be linear: the Pl\"ucker coordinates of congruence lines are related by linear equations. Thus, the classification and geometric description of linearly degenerate 3-component system of Temple class \cite{AF-01} extend to differential-geometric context the algebraic-geometrical classification of linear congruences, obtained by Castelnuovo in 19th century. Although a general classifications of linear congruences in $\mathbb P^5$ is not known, the integrable class of corresponding 4-component systems of conservation laws was described in detail in \cite{AF-05}. These systems are linearly degenerate; their rarefaction curves are rectilinear; their characteristic speeds form a harmonic quadruplet; the corresponding congruence is linear. Integrable systems with 3 or 4 components, corresponding to linear congruences, admit also the third-order nonlocal Hamiltonian formalism \cite{FPV-18}.
It is remarkable that linearly degenerate systems of Temple class exist for any number of component: along those having  linear congruences, there is a parametric family of so-called reducible systems (see \cite{AF-02}) whose congruences are not linear. Such systems are equivalent to completely exceptional equations of Monge-Amp\`ere type \cite{B-92}.

Projective theory of congruences is a rich and fascinating area of classical differential geometry. The above story certifies that its interrelations with mathematical physics are not restricted just to geometrical optics. As a one more example, we mention the interpretation of iterated Laplace transformation of congruences (see \cite{F-50}) as integrable chains for characteristic speeds \cite{AF-99}.

This paper is dedicated to systems (\ref{CL}) admitting at least two additional conservation laws, i.e.
differential one-forms $B(u)dx+A(u)dt$, $N(u)dx+M(u)dt$ closed on each solution of (\ref{CL}). We suppose that
 the $n+4$ one-forms $dx$, $dt$, $u^i(x,t)dx+f^i(u(x,t)dt, \ \ i=1,...,n,$  $B(u)dx+A(u)dt,$ $N(u)dx+M(u)dt$ are linearly independent over $\mathbb R$ (as 1-forms in $(n+2)$-dimensional space with coordinates $u,x,t$).
It is known that generic Hamiltonian systems have exactly two additional conservation laws, and nondiagonalizable systems with 3 or 4 components cannot admit more than two additional conservation laws.

We straightforwardly generalize  the above construction by assigning   to (\ref{CL}) an $n$-dimensional submanifold of Grassmannian $Gr(2,n+4)$:
\begin{equation}
\begin{array}{l}
Y^i=u^i Y^0+f^i(u)Y^{n+3},\ \ \ \ i=1,..., n,\\
\\
Y^{n+1}=B(u) Y^0+A(u)Y^{n+3},\\
\\
Y^{n+2}=N(u) Y^0+M(u)Y^{n+3}.
\end{array}
\label{congrn2}
\end{equation}
Now the lines in $\mathbb P^{n+3}$, corresponding to the points of this submanifold of the Grassmannian $Gr(2,n+4)$, do not sweep out the whole space but form a ruled $(n+1)$-dimensional hypersurface $\Sigma$ in $\mathbb{P}^{n+3}$.
In what follows, we call  $\Sigma$  the {\it ruled hypersurface of conservation law system}. The existence of two additional independent conservation laws is very restrictive and implies non-trivial geometric properties of the hypersurface $\Sigma$. For such systems we show:
\begin{itemize}
\item projective tangent spaces to $\Sigma$ are stable along its line generators;
\item the hypersurface $\Sigma^*$ dual to $\Sigma$ also corresponds do some system of conservation laws, which we call the  {\it dual} system;
\item  Hamiltonian systems of conservation laws are auto-dual;
\item  the ruled hypersurface $\Sigma$ of a Hamiltonian system sits in some quadric, there is a contact structure on  Fano variety
of this quadric, and the generators of $\Sigma$ form a Legendre submanifold;
\item 3-component nondiagonalizable systems of Temple class are dual to systems with constant characteristic speeds.
\end{itemize}
Five conservation laws of 3-component nondiagonalizable systems with constant characteristic speeds can be viewed as the {\it Abelian relations} of 3-webs of curves in the space of field variables $u$, the web being formed by the rarefaction curves. The web rank, i.e. the dimension of the vector space of Abelian relations, is at most 5 for such webs. There is an elegant geometric description of such webs via cubic hypersurfaces in $\mathbb P^4$     (\cite{BB-38,BW-34}). The established duality gives a complete geometric description of 3-component nondiagonalizable systems of Temple class. The construction of 3-component nondiagonalizable systems of Temple class via cubic hypersurfaces in $\mathbb P^4$ was presented in \cite{AF-99}, though it was not proved that this construction embraces all such systems.

We finish the introduction by an example encapsulating practically all the properties mentioned above.
Consider the  Hamiltonian system
$$
\frac{\partial}{\partial t}\left(\begin{array}{l} u_1\\ u_2\\ u_3 \end{array} \right)=\left(\begin{array}{ccc} -1 & 1 & 1\\ 1 & -1 & 1\\ 1 & 1 & -1 \end{array} \right)\frac{d}{dx}\left(\begin{array}{l} \partial_{u_1} H(u)\\  \partial_{u_2} H(u)\\  \partial_{u_3} H(u) \end{array} \right)
$$
with the Hamiltonian density $H(u)=\frac{1}{2}u_1u_2u_3$ and flat coordinates $u_1,u_2,u_3$.

It has 5 conservation laws
$$
\begin{array}{l}
\sigma_1=u_1dx+\frac{1}{2}[u_2u_3-u_1(u_2+u_3)]dt, \\
\\
\sigma_2=u_2dx+\frac{1}{2}[u_3u_1-u_2(u_3+u_1)]dt, \\
\\
\sigma_3=u_3dx+\frac{1}{2}[u_1u_2-u_3(u_1+u_2)]dt, \\
\\
\sigma_4 = b(u)dx+c(u)dt,\\
\\
\sigma_5 = c(u)dx+f(u)dt,\\
\end{array}
$$
where $b=-\frac{1}{2}\left(u_1u_2+u_2u_3+u_3u_1\right)$, $c= u_1u_2u_3$, and $f= b^2-ac$ with $a=u_1+u_2+u_3$.

This system is {\it reducible}, namely, with $\sigma_1+\sigma_2+\sigma_3=a(u)dx+b(u)dt$,  one can introduce a potential $u(x,t)$ (see \cite{MF-96}) by
$$
u_{xxx}=a, \ \ \ u_{xxt}=b, \ \ \ u_{xtt}=c,\ \ \ u_{ttt}=f,
$$
so that the system is equivalent to the associativity equation of 2D topological field theory
$$
u_{ttt}=u^2_{xxt}-u_{xxx}u_{xtt}.
$$

The characteristic speeds of the system are $\lambda_{\alpha}=-u_{\alpha}$, the rarefaction curves corresponding to $\lambda _1=-u_1$ are trajectories of the characteristic vector field $\xi_1=(u_1-u_2)\partial_{u_2}+(u_3-u_1)\partial_{u_3}$. The other two characteristic vector fields are obtained by cyclic permutations of indices. Thus, the system is {\it lineraly degenerate}, i.e. $\xi_{\alpha}(\lambda_{\alpha})=0$. Moreover, in the new field variables $a,b,c$, which are also densities of conservation laws, the rarefaction curves are rectilinear, thus the system is of {\it Temple class} (\cite{AF-96}).

Finally, the reciprocal transformation
$$
dX=\sigma_1-\sigma_3, \ \ \ \ dT=\sigma_2-\sigma_3
$$
brings our system to the one with constant characteristic speeds $\Lambda_1=\infty$, $\Lambda_2=0$, and $\Lambda_3=-1$. Now let us rewrite the system in the new independent variables $X,T$, choose a basis for the conservation laws, and compose two vectors $\vec I_1,\vec I_2$ with 5 components: the first with fluxes (coefficients of $dT$ of corresponding 1-forms) of the elements of the basis and the second with the densities (coefficients of $dX$).  Then the points  $P,Q,R\in \mathbb P^4$, where $P=-\xi_1(\vec I_2)$,
$Q=\xi_2(\vec I_1)$, and $R=\xi_3(\vec I_2)=-\xi_3(\vec I_1)$ parametrize one and the same cubic hypersurface  in $\mathbb P^4$, the three lines spanned by any pair of points  $P,Q,R$ lying on the cubic.

\section{Conservation laws and families of lines in projective space}

Interpretation of  solutions to (\ref{CL})  as surfaces in affine space whose Gauss map images lie in the Grassmannian sub-manifold (\ref{congrn})
 naturally leads to the idea of {\it reciprocal transformation}: instead of $x,t$ as parameters on the surface $S$ (see Introduction) one can choose suitable linear combinations of the affine coordinates $(q_1,...,q_n,x,t)$. Such reparametrizations induce non-local changes of independent variables:
\begin{equation}
\begin{array}{c}
dX=(\alpha _iu^i+\alpha) dx+(\alpha _if^i+\tilde{\alpha}) dt,\\
dT=(\beta _iu^i+\beta) dx+(\beta _if^i+\tilde{\beta} ) dt,
\end{array}
\label{Trec}
\end{equation}
(here $\alpha _i,\alpha,\tilde{\alpha},\beta _i,\beta,\tilde{\beta}$ are arbitrary constants) and the following substitution for the densities $U^i$ and fluxes $F^i$
$$
U^i=\frac{u^iM-f^iN}{BM-AN}, ~~ F^i=\frac{f^iB-u^iA}{BM-AN},
$$
of the transformed system
\begin{equation}
U^i_T=\frac{d}{dX} F^i(U) ,\ \ \,  i=1,...,n,
\label{CLn}
\end{equation}
where $B=\alpha _iu^i+\alpha$, $A=\alpha _if^i+\tilde{\alpha}$, $N=\beta _iu^i+\beta$, amd $M=\beta _if^i+\tilde{\beta}$.

Finally, there is no reason to restrict oneself to affine transformations changing only two of $n+2$ coordinates $(q^1,...,q^n,x,t)$. Thus, the affine transformations of $A^{n+2}$ induce projective transformations of the congruence (\ref{congrn}).

For systems admitting 2 additional conservation laws, we straightforwardly generalize this construction to surfaces  in  the affine space
$A^{n+4}$ with coordinates $(q^1,...,q^n,q^{n+1},q^{n+2},x,t)$, the quantities $q^{n+1},q^{n+2}$ being defined
by $dq^{n+1}=B(u)dx+A(u)dt$ and $dq^{n+2}=N(u)dx+M(u)dt$, and consider all the geometric objects  in $\mathbb P^{n+3}$ up to projective equivalence.

Observe that any $n$-tuple $f^1,...,f^n$ of functions define a system of conservation laws (\ref{CL}), while the existence of two additional independent conservation laws is very restrictive and implies non-trivial geometric properties of the hypersurface $\Sigma$. Following algebraic geometry, we say that $\Sigma$ is {\it nondegenerate} if it does not lie in any hyperplane.

\begin{theorem}\label{tangentstable}
A nondegenerate ruled hypersurface $\Sigma$ of codimension 2 in $\mathbb P^{n+3}$ corresponds via (\ref{congrn2}) to some system of $n$ conservation laws (\ref{CL}), admitting two additional independent conservation laws, if and only if the projective tangent spaces to $\Sigma$ are stable along the line generators (\ref{congrn2}).
\end{theorem}
{\it Proof:} Let $\Sigma$ correspond to a conservation law system.  Line generators of $\Sigma$ are spanned by points $p=[1,u^1:...:u^n:B:N:0]$ and $r=[0:f^1:...:f^n:A:M:1]$. Denoting   the derivations with respect to $u^i$ by subscript $i$, we get the tangent space at $p$ as the span of $p,r,p_i,$ $i=1,...,n$. For the additional conservation laws one has
\begin{equation}\label{AB}
A_i=B_kf^k_i,\ \ \ \ \   M_i=N_kf^k_i,\ \ \ i=1,...,n.
\end{equation}
hence $r_i=f^l_ip_l$. Thus, the projective tangent space at $r$, and therefore at any point of the generator, coincides with the tangent space at $p$.

Conversely, a ruled hypersurface $\Sigma$ of codimension 2 can be locally parametrized by (\ref{congrn2}). Then the stability of tangent spaces implies (\ref{AB}). Therefore $B(u)dx+A(u)dt$ and $N(u)dx+M(u)dt$ are indeed the additional conservation laws of system (\ref{CL}).
\hfill $\Box$\\

A richer geometry corresponds to the {\it strictly hyperbolic} conservation laws systems, for which the matrix $\nabla f=f^i_j=\partial _{u_j}f^i$ has $n$ distinct real eigenvalues $\lambda^{\alpha}$, called {\it characteristic speeds}:
\begin{equation}\label{char}
\nabla f\xi_{\alpha}=\lambda^{\alpha}\xi_{\alpha},
\end{equation}
 (no summation).
The integral curves of the {\it characteristic} vector fields  $\xi_{\alpha}=\xi_{\alpha}^i\partial _i$, called {\it rarefaction curves}, play distinguished role in the theory of hydrodynamic type systems. A generic solution to (\ref{CL}) parametrizes some surface in the space of field variables $u$. When this surface degenerates into a curve then this curve is necessarily tangent to one of $\xi_{\alpha}$.
For the congruence (\ref{congrn}), the rarefaction curves parametrize the {\it developable surfaces} of the congruence: the lines corresponding to points of the rarefaction curve are tangent to the cuspidal edge of the developable surface (see \cite{AF-96}). Applying this correspondence to various projections of the hypersurface  $\Sigma$ onto $n+1$-dimensional projective subspaces, one concludes that rarefaction curves still parametrize developable surfaces, which now lie on  $\Sigma$. Thus, each generator of  $\Sigma$ is tangent to $n$ submanifolds of  $\Sigma$, which we will call {\it focal}. These focal submanifolds are parametrized by $u$ as follows (see again see \cite{AF-96}):
\begin{equation}
\begin{array}{l}
Y^0=-\lambda^{\alpha}Y^{n+3}\\
\\
Y^i=(f^i(u)-\lambda^{\alpha} u^i)Y^{n+3},\ \ \ \ i=1,..., n,\\
\\
Y^{n+1}=(A(u)-\lambda^{\alpha}B(u))Y^{n+3},\\
\\
Y^{n+2}=(M(u)-\lambda^{\alpha}N(u))Y^{n+3}.
\end{array}
\label{focal}
\end{equation}

\section{Dual conservation laws}

Consider the hypersurface dual to $\Sigma$. Recall that the dual hypersurface $\Sigma ^*$ is the set of points $h$ in the dual space $\mathbb (\mathbb P^{n+3})^*$ such that each $h$ is a projective tangent hyperplane to $\Sigma$.
Note that $codim(\Sigma)=2$ and therefore for each $p\in \Sigma$, there is a pencil $L$ of projective  hyperplanes tangent to $\Sigma$ at $p$.
Thus, the dual hypersurface $\Sigma ^*$ is also ruled.
Theorem \ref{tangentstable} implies that  the pencils $L$ are stable along the line generators of $\Sigma$.
Consequently, the family of the line generators $L$ of the dual hypersurface $\Sigma ^*$ can be parametrized by $u$ and therefore the hypersurface $\Sigma ^*$ has the codimension at least 2.
Observe also that n+2 one-forms of the conservation laws are independent if and only if the hypersurface is not contained in a proper projective subspace of $\mathbb{P}^{n+3}$.
\begin{theorem}\label{dual}
Suppose that the ruled hypersurface $\Sigma$, corresponding via (\ref{congrn2}) to some system of $n$ conservation laws (\ref{CL}), admitting two additional independent conservation laws, enjoys the following properties:\\
1) it is not ruled by projective spaces of any dimension $r>1$,\\
2) its projective tangent spaces are not concurrent.\\
Then the dual hypersurface $\Sigma ^*$ is of codimension 2, it is not contained in a proper projective subspace of $\mathbb{P}^{n+3}$, and its projective tangent spaces are stable along the generators of $\Sigma ^*$.
\end{theorem}
{\it Proof:} By the biduality theorem, the dual to $\Sigma ^*$ is $\Sigma $. If $\Sigma ^*$ were of codimension larger than $2$ then $\Sigma $ would allow ruling in projective spaces of dimension at least 2 (see e.g. \cite{GKZ-94}). If  $\Sigma ^*$ were contained in a proper projective subspace of $\mathbb{P}^{n+3}$ then the tangent spaces of
$\Sigma $ were concurrent.
To prove the last claim, we introduce homogeneous coordinates $Z_i,\ i=0,...,n+3$ in the dual space so that the incidence is given by $Z_iY^i=0.$ Direct calculation gives the following ruling of $\Sigma ^*$:
\begin{equation}
\begin{array}{l}

Z_0=(u^iB_i(u)-B(u))Z_{n+1}+(u^iN_i(u)-N(u))Z_{n+2},\\
\\
Z_i=-B_i(u) Z_{n+1}-N_i(u)Z_{n+2},\ \ \ \ i=1,..., n,\\
\\
Z_{n+3}=(f^i(u)B_i(u)-A(u))Z_{n+1}+(f^i(u)N_i(u)-M(u))Z_{n+2}.\\
\end{array}
\label{congrn2D}
\end{equation}
The generators are spanned by
$\tilde{p}=[u^iB_i-B:-B_1:...:-B_n:1:0:f^iB_i-A]$ and $\tilde{r}=[u^iN_i-N:-N_1:...:-N_n:0:1:f^iN_i-N]$.

Observe that $\partial_k(u^iB_i-B)=u^iB_{ik}$ and by (\ref{AB}) holds  $\partial_k(f^iB_i-A)=f^iB_{ik}.$
Similarly, $\partial_k(u^iN_i-N)=u^iN_{ik}$ and  $\partial_k(f^iN_i-N)=f^iN_{ik}.$ Thus, under the  restrictions imposed, the hypersurface $\Sigma ^*$
  also comes from some conservation law system, for which either $B_i$ or $N_i$ play the role of field variables (densities). Hence the projective tangent space to  $\Sigma ^*$ along the line $\tilde{p}\tilde{r}$ is spanned by $\tilde{p},\tilde{r},\tilde{p_i},\ i=1,...,n$ and is stable.
 \hfill $\Box$

\begin{definition}
Two systems of conservation laws, admitting 2 additional conservation laws, are dual to each other if their ruled hypersurfaces are dual.
\end{definition}

\noindent{\bf Remark.}
Exploiting the geometric ideas manifested in Theorems \ref{tangentstable} and \ref{dual}, we easily construct  systems of conservation laws, which admit one additional conservation law,   via envelope of  an $n$-parameter family of hyperplanes in $\mathbb P^{n+2}$
$$
Z_0+u^iZ_i+B(u)Z_{n+1}+N(u)Z_{n+2}=0,\ \ \ i=1,...,n.
$$
The corresponding $n$-parameter family of lines is
$$
\begin{array}{l}
Z_0=(u^iB_i-B)Z_{n+1}+(u^iN_i-N)Z_{n+2},   \\
\\
 Z_i=-B_i Z_{n+1}-N_iZ_{n+2},\ \ \ \ i=1,..., n.\\
\end{array}
$$
The constructed systems are {\it Godunov's systems} for the case of 2 independent variables  \cite{Gi-61}.
\medskip

 For a strictly hyperbolic system of hydrodynamic type, one can introduce {\it characteristic  forms} $\omega^{\alpha}=l_i^{\alpha}du^i$, where $l^{\alpha}(u)=(l^{\alpha}_1(\vec{u}),...,l^{\alpha}_n(\vec{u}))$ are left eigenvectors of the matrix $\nabla f$. Then the system (\ref{CL}) (and, similarly, any strictly hyperbolic system of hydrodynamical type \cite{F-93}) can be expressed as $n$ exterior equations (no summation)
\begin{equation}\label{exterior}
\omega^{\alpha}\wedge (dx+\lambda^{\alpha}dt)=0.
\end{equation}

This form is particularly convenient for studying systems without Riemann invariants. The forms $\omega^{\alpha}$ constitute a basis for each cotangent space of dependent variables. Thus, differentials of functions $A,B,M,N$ of dependent variables can be decomposed as $dA=A_{\alpha}\omega^{\alpha}$ etc.
Then the equations (\ref{AB}) for conservation laws  become (no summation)
\begin{equation}\label{ABf}
A_{\alpha}=\lambda^{\alpha} B_{\alpha},\ \ \ \ \   M_{\alpha}=\lambda^{\alpha} N_{\alpha},\ \ \alpha=1,...,n.
\end{equation}
 If we normalize the characteristic forms and the characteristic vector fields to verify $\omega^{\alpha}(\xi_{\beta})=\delta^{\alpha}_{\beta}$ then $A_{\alpha}=\xi_{\alpha}(A)$.
In particular, we define coefficients $c^{\alpha}_{ij}$ and  $\lambda^{\alpha}_i$ by relations
 \begin{equation}\label{structureeq}
d\omega ^{\alpha}=c^{\alpha}_{ij}\omega ^i\wedge \omega ^j, \ \ \ \ d\lambda^{\alpha}=\lambda^{\alpha}_i\omega^i.
\end{equation}
and call the first ones  {\it structure equations}.

\begin{theorem}\label{structureD}
Characteristic vectors $\xi_{\alpha}$ of strictly hyperbolic system (\ref{CL}) are also characteristic vectors of the dual system. In particular, the structure equations of the dual system are the same.
\end{theorem}
 {\it Proof:} Suppose that $B_i$ are functionally independent densities of the the dual system defined by (\ref{congrn2D}).  Differentiating the first equation of (\ref{AB}) and comparing the mixed derivatives, one gets
 $$
 A_{ij}=B_{kj}f^k_i+B_kf^k_{ij}=A_{ji}=B_{ki}f^k_j+B_kf^k_{ij}.
 $$
Therefore $B_{kj}f^k_i=B_{ki}f^k_j$ and for each characteristic vector  $\xi_{\alpha}$ of (\ref{CL}) holds
$$
\xi^i_{\alpha}B_{ki}f^k_j=B_{kj}f^k_i\xi^i_{\alpha}=B_{kj}\lambda^{\alpha}\xi^k_{\alpha}.
$$
Note that $\xi^i_{\alpha}B_{ki}=\xi_{\alpha}(B_{k})$ are the components of $\xi_{\alpha}$ in coordinates $B_k$.
Then the above equality says that $\xi_{\alpha}(B_{k})$ are coordinates of a (left) eigenvector of $\nabla f$. Repeating the calculation for $N$ and taking into account the strict hyperbolicity of (\ref{CL}) we conclude the proof.  \hfill $\Box$\\
\section{Contact structure on Fano variety of quadric and Hamiltonian systems}

A Hamiltonian system of hydrodynamical type of first order (see \cite{DN-83} for details) can be written in the so-called flat coordinates as
\begin{equation}\label{hamiltonian}
u^i_t=\varepsilon^i\frac{d}{dx}\left(\frac{\partial h(u)}{\partial u^i}\right),\ \ \  i=1,...,n,
\end{equation}
(no summation) where $h(u)$ is  the Hamiltonian density and $\varepsilon^i=\pm 1$. For $n\ge 3$, a generic system (\ref{hamiltonian}) has no Riemann invariants and possesses exactly $n+2$ conservation laws \cite{T-85}:
\begin{equation}
\begin{array}{l}
u^idx+\varepsilon^i h_idt, \ \text{(no summation)}\\
\\
\frac{1}{2}|u|^2dx+(h_iu^i-h)dt,\\
\\
hdx+\frac{1}{2}|\nabla h|^2dt,
\end{array}
\end{equation}
where $h_i=\frac{\partial h}{\partial u^i}$, $\nabla h=(h_1,...,h_i)$, and $|v|=\varepsilon^i(v_i)^2$ is a pseudo-Riemannian metric.

The line generators of the ruled hypersurface $\Sigma$ are spanned by points
$$
p=[1:u^1:...:u^n:|u|^2/2:h:0]\ \mbox{and} \
r=[0:\varepsilon^1 h_1,...,\varepsilon^n h_n:h_iu^i-h:|\nabla h|^2/2:1].
$$
It is straightforward that these lines lie on the smooth quadric  $(X,X)=0$ defined by the following symmetric bilinear form
$$
(X,Y)=\varepsilon^1X_1Y_1+...+\varepsilon^nX_nY_n-X_0Y_{n+1}-X_{n+1}Y_0-X_{n+2}Y_{n+3}-X_{n+3}Y_{n+2}.
$$
The lines $pr$ parameterize an $n$-dimensional submanifold of $(2n+1)$-dimensional Fano variety $F_{1,n+2}$ of lines on this $(n+2)$-dimensional quadric  $Q$. The tangent space to $F_{1,n+2}$ at $l\in F_{1,n+2}$, where $l$ is the projectivization of a plane $\Lambda \subset \mathbb R^{n+4}$, can be represented as
$$
T_lF_{1,n+2}=\{\overline{\varphi}\in \mbox{Hom}(\Lambda, \mathbb R^{n+4}/ \Lambda): \overline{\varphi}(v)\in T_{[v]}Q/ \Lambda) \},
$$
 (see \cite{H-92},  p. 209). For a representative $\varphi\in \mbox{Hom}(\Lambda, \mathbb R^{n+4})$ of $\overline{\varphi}$ and $X,Y\in \Lambda$, the imposed restrictions read as
$$
(X,\varphi(X))=(Y,\varphi(Y))=(X,\varphi(Y))+(Y,\varphi(X))=0.
$$
Thus, the additional constraint
$$
(X,\varphi(Y))=0
$$
defines a distribution $\tau$ of codimension one on $TF_{1,n+2}$. The geometrical meaning of this constraint is transparent: any pair of points $[X],[Y]$ on a line generator of the quadric $Q$ determines a pencil of projective tangent hyperplanes to $Q$ along the generator, and    $\varphi$ sends any point $Z\in \Lambda$ to the codimension 2 plane of intersection of pencil hyperplanes.

For $\varepsilon^i=1,\ \ i=1,...,n$ and up to a linear change of variables $Y_0,Y_{n+1},Y_{n+2},Y_{n+3}$, the quadric $Q$ coincides with the {\it Lie quadric} in Lie sphere geometry, the distribution $\tau$ defining the contact structure on the Fano variety of Lie quadric (see \cite{C-08}). The introduced distribution is a straightforward generalization.
\begin{lemma}
The distribution $\tau$ provides $TF_{1,n+2}$ with a contact structure.
\end{lemma}
Observe that with the points $p,s$
$$
p=[1:u^1:...:u^n:|u|^2/2:h:0]\ \mbox{and} \
s=[0:\varepsilon^1 k_1,...,\varepsilon^n k_n:k_iu^i-h:|\nabla h|^2/2:1],
$$
and $(u,k,h)\in \mathbb R^{2n+1}$, the lines $ps$ parameterize a Zarisski open subset of $F_{1,n+2}$.

\begin{theorem}
A ruled hypersurface corresponds to a Hamiltonian system (\ref{hamiltonian}) if and only if its line generators form a Legendre submanifold of Fano variety $F_{1,n+2}$ of some smooth quadric in $\mathbb P^{n+3}$.
\end{theorem}
{\it Proof:} The claim follows from the calculation: $(s,dp)=k_idu^i-dh$. \hfill $\Box$\\ \smallskip

 For Hamiltonian systems, duality does not produce any new class of systems with 2 additional conservation laws.
\begin{theorem}\label{autodual}
Hamiltonian system of conservation laws (\ref{hamiltonian}) is auto-dual, points of $\Sigma^*$ being the polar hyperplanes of the corresponding points of $\Sigma$ with respect to the quadric fixed by (\ref{hamiltonian}).
\end{theorem}
{\it Proof:} Applying formulae  (\ref{congrn2D}) to $f^i=\varepsilon^i h_i$ (no summation), $A=h_iu^i-h$, $B=\frac{1}{2}|u|^2$, $M=\frac{1}{2}|\nabla h|^2$, $N=h$, and substituting $Z_{n+3}=Y_{n+2}$, $Z_{n+1}=Y_0$, $Z_{n+2}=Y_{n+3}$, $Z_0=Y_{n+1}$, $ Z_i=-\varepsilon^i Y_i$  (no summation), one arrives at (\ref{congrn2}).
\hfill $\Box$\\

\section{Maximal rank 3-webs of curves in space and conservation laws}\label{5}

Consider 3 foliation $\mathfrak{F}_1,\mathfrak{F}_2,\mathfrak{F}_3$ of space by curves. They form a  {\it 3-web} $\mathfrak{W}_3=\{\mathfrak{F}_1,\mathfrak{F}_2,\mathfrak{F}_3 \}$. A function  $I_k$ (defined, generically, only locally) is a {\it first integral} of the $k$th foliation $\mathfrak{F}_k$ of the web, if it is constant on the leaves of $\mathfrak{F}_k$. An Abelian relation of $\mathfrak{W}_3$ is  an identity relating the first integrals of the form
\begin{equation}\label{Abel}
dI_1+dI_2+dI_3\equiv0.
\end{equation}

For a given web, Abelian relations form a vector space, whose dimension ${\rm rk}(\mathfrak{W}_3)$ is called the web rank. It is known (see \cite{BB-38}) that a $3$-web of integral curves of  3 vector fields have the rank at most 5 if the 2-dimensional distributions spanned by vector field pairs are not integrable. There is the following elegant description of such 3-webs of maximal rank (see \cite{BB-38} p.313 or \cite{BW-34}) via  cubic hypersurfaces $M^3$ in projective space $\mathbb P^4.$ A generic 2-dimensional plane $\mathbb P^2\subset \mathbb P^4$ cuts $M^3$ along some cubic curve $C$.
 There is a 3-parameter family $F^3$ of planes for which the curve  $C$ splits into 3 lines $l_1,l_2,l_3$.  Each line $l_{\alpha}$ determines a one parameter family $F^1_{\alpha}\subset F^3$ of planes, all planes  of $F^1_{\alpha}$ containing  the line $l_{\alpha}$.  These 3 families $F^1_{\alpha},\ \ \alpha=1,2,3,$ are the foliations of some 3-web. This 3-web is of maximal rank and each maximal rank 3-web of curves in space is locally diffeomorphic to some web obtained by the above construction. Observe that the projective moduli space of cubic hypersurfaces $M^3\subset \mathbb P^4$ is of dimension 10, thus the 3-webs under study depend on 10 essential parameters.

 Now let us choose some local parameters  $u^1,u^2,u^3$ in the family $F^3$. Then at each point $u=(u^1,u^2,u^3)$ in the parameter space we can (also locally) define 3 vector fields $\xi_{\alpha}$, each tangent to the corresponding family $F^1_{\alpha}$.
  Integrating an Abelian relation (\ref{Abel}), one can choose the integration constants for $I_i$ to hold $I_1+I_2+I_3\equiv0$. We have $\xi_1(I_1)=0$, $\xi_2(I_2)=0$,  and $\xi_3(I_3)=-\xi_3(I_1)-\xi_3(I_2)=0$.
These three relations can be interpreted as equations (\ref{ABf}) for the conservation laws $I_1dt+I_2dx$ of the hydrodynamical type system (\ref{exterior}) with characteristic speeds $\lambda_1=\infty$, $\lambda_2=0$, and $\lambda_3=-1$ (if you feel uncomfortable with infinity, rotate the coordinate axis in $xt$-plane).
The system has 5-dimensional space of conservation laws.
\begin{lemma}\label{PQR}
 The ruled hypersurface of the above system of conservation laws is 4-dimensional.
\end{lemma}
 {\it Proof:} The proof appeals to the description of 3-webs of maximal rank in \cite{BB-38}. Namely, we choose some basis of Abelian relations $(I_1^i,I_2^i,I_3^i)$, $i=1,...,5$. Then  the 3 lines $l_1,l_2,l_3$  generically intersect in 3 points $P,Q,R\in \mathbb P^4$, where $P=\xi_1(\vec I_3)=-\xi_1(\vec I_2)$,
$Q=\xi_2(\vec I_1)=-\xi_2(\vec I_3)$, and $R=\xi_3(\vec I_2)=-\xi_3(\vec I_1)$. Rotating the axis in the $xt$-plane, we obtain a system of conservation laws with densities $I_1^i-I_2^i$. Then $\xi_1(\vec I_1- \vec I_2)=P$, $\xi_2(\vec I_1-\vec I_2)=Q$, $\xi_3(\vec I_1-\vec I_2)=-2R$ are linearly independent. Thus, the densities of some 3 conservation laws of the transformed system are functionally independent. Therefore the ruled hypersurface of the transformed system, as well as of the original one, is 4-dimensional.   \hfill $\Box$\\

The above Lemma implies that the  systems of conservation laws corresponding to maximal rank 3-webs are not degenerate. The characteristic speeds are not invariant under reciprocal transformations, though one can geometrically express  the property of being equivalent to a system with constant characteristic speeds   as follows.
\begin{theorem}
Strictly hyperbolic system of conservation laws is  equivalent by reciprocal transformation to one with constant characteristic speeds if and only if its focal manifolds lie in hyperplanes.
\end{theorem}
 {\it Proof:} The hyperplane equation is the first one in (\ref{focal}). \hfill $\Box$\\

 We can reformulate the result of Blaschke and Walberer \cite{BW-34} as follows:
\begin{theorem}
If the focal manifolds of strictly hyperbolic system of conservation laws with 2 additional conservation laws  lie in hyperplanes then
the system can be obtained from a maximal rank 3-web of curves in space as described above.
\end{theorem}

\section{Systems of Temple class}

Generic solutions to hyperbolic systems of hydrodynamic type develop shock waves (see \cite{L-79}). For a system of conservation laws (\ref{CL}),  the values of field variables $u_0$ before the shock wave, which is a curve in $xt$-plane, and their values $u$ across the shock wave are related by the {\it shock curve} (in $u$-space) with the vertex at $u_0$. To the shock curve with the vertex $u_0$, there corresponds a 2-dimensional ruled surface of the congruence (\ref{congrn}). The generators of this surface are the congruence lines $l(u)$ intersecting the "vertex" line $l(u_0)$ (see \cite{AF-96}.
\medskip

\noindent{\bf Remark.}
Note that conditions on the shock wave, being obtained by integrating the one-forms $u^idx+f^idt$ over an infinitesimal contour along the arc of the shock wave, distinguish the n-dimensional subspace of conservation laws, spanned by those with densities $u^i$. The choice of quantities preserved also in "integral" form, is dictated by the physics of the model in question. For systems with additional conservation laws,  these additional quantities are preserved only for smooth solutions.
\medskip

As shown by Lax \cite{L-57}, a shock curve with the vertex at a generic point $u_0$ splits into $n$ branches, the $\alpha$-th branch being
$C^2-$tangent of some rarefaction curve of the $\alpha$-th family passing through $u_0.$
There are also systems with shock curves coinciding  with their associated
rarefaction curves. For such systems, among smooth solutions with degenerate (i.e. one-dimensional) hodograph, there exit also weak, discontinuous solutions whose hodograph belongs to a smooth curve.    Systems with coinciding shock and rarefaction curves were characterized by Temple \cite{T-83} as follows.
\begin{theorem}\cite{T-83}
Rarefaction curves of the $\alpha$-th family coincide
with the associated  branches of the shock curve if and only if either\\
1) every rarefaction curve of the $\alpha$-th family is a straight line in the $u-$space, or\\
2) the characteristic velocity $\lambda ^{\alpha}$ is constant along rarefaction curves of the $\alpha$-th family:
 $$
 \xi_{\alpha}(\lambda ^{\alpha})=0
 $$
 (no summation).
\end{theorem}
Systems enjoying the second condition, are called {\it linearly degenerate}. Such systems possess properties unusual for nonlinear hyperbolic systems, namely, their smooth solutions do not develop shock waves \cite{R-67}.  Reciprocal transformations do not preserve the characteristic speeds, but respect linear degeneracy  \cite{F-89}.
Geometrically, the condition means that
developable surfaces of the congruence (\ref{congrn}) (see \cite{AF-96}), or of the hypersurface (\ref{congrn2}) respectively, are conical, and therefore the codimension of  focal submanifolds jumps up at least by one.

The first condition of Temple's Theorem also translates in the geometrical language quite naturally: the cuspidal edges of the developable surfaces are planar \cite{AF-96}.
This geometrical condition is respected not by all reciprocal transformations, but by ones corresponding to  linear combinations of "canonical" conservation laws $u^idx+f^idt$ (see also \cite{AF-96}).
In what follows, we call systems of conservation laws with rectilinear rarefaction curves {\it systems of Temple class}.

In this section, we will be concerned with {\it nondiagonalizable} 3-component systems of Temple class not admitting a single Riemann invariant, which amounts to $d\omega ^{\alpha} \wedge \omega ^{\alpha}\ne 0$ (no summation) for $\alpha =1,2,3$. It is well known that nondiagonalizable 3-component systems admit at most 5 conservation laws.
\begin{theorem}\label{5temple}
Nondiagonalizable 3-component systems of Temple class admit 5 conservation laws.
\end{theorem}
This Theorem can be proved on the basis of computations made in \cite{AF-99}: it was shown that the condition of having rectilinear rarefaction curves is very rigid, the imposed restriction almost  fixes (up to some constant) the coefficients $c^{\alpha}_{ij}$ and $\lambda^{\alpha}_i$ of equations (\ref{structureeq}).
In fact, these equations allow to compute the dimension of the space of conservation laws.
We can also get Theorem \ref{5temple} as a corollary of the proposed duality.
\begin{theorem}\label{DT3}
Nondiagonalizable 3-component systems of Temple class are dual to the systems  obtained from maximal rank 3-webs of curves in space.
\end{theorem}
 {\it Proof:} Let $(I_1^i,I_2^i,I_3^i)$, $i=1,...,5$ be a basis of Abelian relations. Then the ruled hypersurface of corresponding system is parametrized by $Y^i=I_1^iY^0+I_2^iY^6$, $i=1,...,5$. Denote $U=(1,I_1^1,...,I_1^5,0)$, $V=(0,I_2^1,...,I_2^5,1).$ With $P,Q,R$ introduced in the proof of Lemma \ref{PQR}, coordinates $Z_i$, $i=0,1,...,6$ in the dual space, and $\vec z=(Z_1,...,Z_5)$, we get the following equations for generators of the ruled hypersurface for the dual system:
$$
\begin{array}{l}
Z_0=-\vec I_1\cdot \vec z,\\
\\
Z_6=-\vec I_2\cdot \vec z,\\
\\
P\cdot \vec z=Q\cdot \vec z=R\cdot \vec z=0.
\end{array}
$$
Projecting this hypersurface to the projective subspace $Z_0=Z_6=0$, we get a line congruence in $\mathbb P^4$ whose   lines are dual to 2-dimensional planes spanned by 3 points $P,Q,R\in \mathbb P^4$. To the rarefaction curves, there correspond the developable surfaces, whose line generators are dual to the planes in one-parameter families denoted by $F^1_{\alpha}$ in section \ref{5}. The planes in the families pass through a pair of points $P,Q,R$, this pair being fixed by the rarefaction curve. This means that the developable surfaces of the congruence are planar and the dual system is of Temple class.
    \hfill $\Box$\\

Invoking the description of Blaschke and Walberer, we obtain (at least locally)  the rectilinear rarefaction curves of the system of Temple class as follows: consider the planes spanned by 3 points $P,Q,R$; take the dual objets, which are lines in $\mathbb P^4$,  and project these lines from a point to a subspace $\mathbb P^3$. The line images are rectilinear rarefaction curves.  This construction of systems of Temple class was proposed in \cite{AF-99}. Theorem \ref{DT3} implies that the construction gives {\it all} nondiagonalizable 3-component  systems of Temple class.

\section{Concluding remarks}

\subsection{Hydrodynamic surfaces} Formulae (\ref{congrn2D}) suggest that we can use the same field variables $u$ to parametrize the ruled hypersurfaces for both the original system and its dual. Due to Theorem \ref{autodual},  duality respects hydrodynamic surfaces of Hamiltonian systems. A general system with 2 additional conservation laws shares the rarefaction curves with its dual (see Theorem \ref{structureD}), thus respecting degenerate hodographs. For a non-Hamiltonian system, its hydrodynamic surfaces are not the ones for the dual system, as  can be shown by calculation for systems of Temple class with the help of equations (\ref{structureeq}).

\subsection{Intersection of classes of described examples}
We have described three sets of 3-component systems with 5 conservation laws: nondiagonalizable Hamiltonian, nondiagonalizable of Temple class, and the systems obtained from maximal rank 3-webs of curves in space. Up to reciprocal transformations, there is only one systems in the intersection of all three sets, namely, the one equivalent to the associativity equation. In fact, the reciprocal transformation orbit of this system is also the intersection of any pair of the described classes. For the 1st and 2nd, it follows from results of \cite{AF-99}: a nondiagonalizable system of Temple class is Hamiltonian if and only if it is linearly degenerate, and, up to reciprocal transformation, there is only one nondiagonalizable Hamiltonian system of Temple class. For the 2nd and 3d, the intesection consists of linearly degenerate systems and we are done by the same reason. Finally, the systems in the intersection of the 3d and the 1st sets are again linearly degenerate and the claim is equivalent to the main Theorem of \cite{F-95}.

\subsection{Reducible systems, linear congruences, nonlocal Hamiltonian}
Linearly degenerate 3-component systems of Temple class can be realized as {\it reducible} systems (see \cite{A-98,AF-02}). The corresponding congruences (\ref{congrn}) turned out to be linear (see \cite{AF-01}). The congruence of the following 4-component system
$$
u^1_t=u^2_x,\ \ \  u^2_t=u^3_x,\ \ \ \ u^3_t=u^4_x,\ \ \ \ u^4_t=\frac{d}{dx}[(u^2)^2-u^1u^3]
$$
is linear, the system is reducible, but it does not admit any additional conservation laws. Moreover, this system possesses nonlocal Hamiltonian structure of third order \cite{FPV-18}. Thus, the known classes of 3-component systems with 2 additional conservation laws are not straightforwardly generalizable.

\subsection{Correlations and dual congruences in $\mathbb P^3$}
Two-component systems of conservation laws have infinitely many additional conservation laws, but duality for such systems can be defined directly from the congruence (\ref{congrn}). Since the object dual to a line in $\mathbb P^3$ is again a line, any congruence has the dual one. If a developable surface of a congruence is conical, the corresponding developable surface of the dual congruence is planar. This kind of duality for 2-component systems, which is different from the one proposed in this paper, was described in \cite{AF-96} via {\it correlations} of $\mathbb P^3$.

\section*{Acknowledgements}
The author thanks E.V. Ferapontov for useful discussions.
This research was supported by FAPESP grant \#2018/20009-6.

\end{document}